# ANALYSIS OF FIRST ORDER SYSTEMS FOR THE SOLUTION OF LAPLACE'S EQUATION

*Vance Faber*
*April 13, 2014*

**Abstract.** *In [3], a randomized iterative method is given that produces an algorithm that has expected running time nearly linear (that is, time less than $O(m^{1+\varepsilon})$ where $m$ is the size of the problem) for approximating the solution to the discrete Laplace equation. The goal of this paper is to explain this algorithm in the language of difference operators on graphs.*

**Introduction**. We analyze a method for solving discretizations of Laplace's equation by splitting into two first order difference equations. This method was first considered in [1] as a way to achieve second order accuracy on irregular grids. In [2], we discussed the importance to this method of functions that sum to zero on cycles. In [3], a randomized iterative method is given that produces an algorithm that has expected running time nearly linear (that is, time less than $O(m^{1+\varepsilon})$ where $m$ is the size of the problem) for approximating the solution to the discrete Laplace equation by splitting in precisely this way. The goal of this paper is to explain this algorithm in the language of difference operators on graphs.

**Laplacian on graphs**. We start with an arbitrary connected directed graph, $G = (V, E)$ with no loops or cycles of length 2. In this case, $E$ is a set of pairs of distinct vertices in $V$ and only one of $(a,b)$ or $(b,a)$ is an edge. Given any real function $u$ on $V$, we define the *gradient* of $u$, $\nabla u$, as the function on the edges given by

$$\nabla u(a,b) = u(b) - u(a).$$

The gradient is a linear operator. The transpose is the *divergence* operator that maps functions on edges to functions on vertices. If $A$ is a function on $E$ then

$$\nabla \cdot A(a) = \sum_{(b,a) \in E} A(b,a) - \sum_{(a,b) \in E} A(a,b).$$

That is, the divergence sums the function on the edges around a given vertex.

The *Laplacian* on the graph is then defined to be the symmetric operator mapping functions on $V$ to functions on $V$ defined by

$$\nabla \cdot \nabla u(a) = \sum_{(a,b) \in E} (u(a) - u(b)) + \sum_{(b,a) \in E} (u(a) - u(b)).$$

Note that the Laplacian is given simply as the difference between the diagonal matrix of degrees of the vertices and the adjacency matrix of the symmetric version of $G$.

**Relationship to discretizations.** If we start with the discretization of Laplace's equation on a mesh, we obtain the discrete Laplacian on the graph that represents the mesh except for additional scale factors related to distances between the vertices. We suppress these factors. We also completely ignore boundary conditions. These additional considerations complicate the operators but do not substantially effect the conclusions.

**The curl and the space of divergence free functions.** We start with the Laplacian equation

$$\nabla \cdot \nabla u = f.  \qquad \text{(Equation 1)}$$

Given any two functions $f$ and $g$ on the same finite domain $D$, we let $(f,g)_D$ be the inner product, the sum of the product of the corresponding entries. We suppress the domain if it is known by context. As usual, we let $\|f\|^2 = (f,f)$. Let $\vec{1}$ be the function that is unity on every vertex of $G$. Note that if Equation 1 has a solution then

$$(\vec{1}, f) = (\vec{1}, \nabla \cdot \nabla u) = (\nabla \vec{1}, \nabla u) = 0.$$

We assume throughout that $(\vec{1}, f) = 0$.

Given $f$, we want to solve Equation 1 for $u$ by utilizing the first order system

$$\begin{aligned} \nabla \cdot A &= f \\ A &= \nabla u. \end{aligned} \qquad \text{(Equation 2)}$$

The method we use will to be to find a function that satisfies the divergence equation (this is easy; see below) and then add to it a function that is divergence free to satisfy the gradient equation. It is classical result that a vector valued function in $\Re^n$ is a gradient if and only if it sums to zero around any closed loop. We use the analogous graph theoretic property (see Lemma 1 below) to find the solution. To this end, we define the *curl* operator, $\nabla \times A$, on functions on the edges. The curl is a dual operator that maps any function $A$ on the edges to a function on the divergence free functions.

The divergence free functions form a finite dimensional vector space. We will show later that the dimension of this subspace of all functions on the edges is $p = m - n + 1$ where $n$ is the number of vertices and $m$ is the number of edges. If $\Gamma$ is a cycle in the undirected graph associated with $G$, suppose that we write $\Gamma$ as $(a_1, a_2, \cdots, a_k, a_1)$. Then the function $C = C_\Gamma$ given by (the indices are read modulo $k$)

$$C(a_i, a_{i+1}) = 1 \quad \text{if } (a_i, a_{i+1}) \text{ is an edge in } E$$
$$C(a_{i+1}, a_i) = -1 \text{ if } (a_{i+1}, a_i) \text{ is an edge in } E$$
$$C(e) = 0 \quad \quad \text{for every other edge}$$

is divergence free. Similarly, we can define a divergence free function with values plus or minus 1 if $\Gamma$ is a disjoint union of cycles. We can choose a basis for the divergence free functions that consist entirely of functions of the type $C_\Gamma$. We call this basis a cycle basis. (See [4] for more complete information about cycle bases.) Then to define the *curl*, it is only necessary to define it on basis functions. We define

$$\nabla \times A(C_\Gamma) = (A, C_\Gamma).$$

*Definition.* Given a function $A$ on the edges of a tree $T$ with each edge having a unique direction and a leaf $s$, we define a function $v = v_{T,s}$ on the vertices of the tree as follows. Let $v(s) = 0$. Then for any vertex $t$, let $P$ be a path from $s$ to $t$ in the tree. We define $u(t)$ to be the sum of $A$ along the edges of the path with a plus sign if the edge is directed toward $t$ and a minus sign if it is not. We say this function is *induced* from $A$ by $T$ and $s$.

*Lemma* 1. *The curl of $A$ is zero if and only if $A$ is a gradient.*

*Proof.* Let $A = \nabla u$ and $C$ be a cycle in the cycle basis. Then $\nabla \times A(C) = (\nabla u, C)$ which is clearly zero since for each vertex $a$ on the cycle $u(a)$ appears twice, once with a plus sign and once with a minus sign.
    Conversely, suppose $\nabla \times A = 0$. Then let $T$ be any spanning tree of $G$ and $s$ be any leaf. If $u(t)$ is the function induced from $A$ by $T$ and $s$, The value $u(t)$ defined in this way is independent of the tree because if $Q$ is a second path from $s$ to $t$ then the path $C$ from $s$ to $t$ along $P$ followed by the path from $t$ to $s$ along $Q$ has even degree at each vertex and thus is in the cycle space of the graph. Since $0 = \nabla \times A(C) = (A, C)$, it must be that the sum of $A$ along the edges of the path $P$ is equal to the sum of $A$ along the edges of the path $Q$. Now let $(a, b)$ be any directed edge in the graph and let $P$ be a path from $s$ to

$a$. Then by considering the path to $b$ along $P$ and then along $(a,b)$, we see that $u(b) = u(a) + A(a,b)$. In other words, $A(a,b) = u(b) - u(a) = \nabla u(a,b)$.

We can now write the system in Equation 2 as

$$\nabla \cdot A = f \qquad \text{(Equation 3)}$$
$$\nabla \times A = 0.$$

*Lemma 2. Consider the extremal problem given by*

$$\text{minimize } \|A\|^2 \text{ subject to } \nabla \cdot A = f. \qquad \text{(Equation 4)}$$

*This has the same solutions as Equation 3.*

Proof. First, we show that the equation has a feasible solution $A_f$. This solution is based on a spanning tree $T$. The function $A_f$ will be zero off the edges of $T$ and satisfy $\nabla \cdot A_f = f$ on the edges (and thus in the whole graph.) We show the existence of $A_f$ recursively on arbitrary trees. To start, consider a tree with just one edge $(a,b)$. Then in order for there to be a solution $0 = (\vec{1}, f) = f(a) + f(b)$. We define $A_f(a,b) = f(b)$. Then the divergence definition gives $\nabla \cdot A_f(a) = -A_f(a,b) = -f(b) = f(a)$ and $\nabla \cdot A_f(b) = A_f(a,b) = f(b)$ as required. Now consider an arbitrary edge at a leaf node $a$. Then either $(a,b)$ or $(b,a)$ is an edge. Suppose that $(a,b)$ is an edge. Then we must have $A_f(a,b) = -f(a)$. (Similarly, if $(b,a)$ is an edge, then $A_f(a,b) = f(a)$.) We now create a function $g$ on a smaller tree $T^- = T - (a,b)$. The function $g$ agrees with $f$ everywhere except $g(b) = f(b) + f(a)$. Note that $(\vec{1}, g) = (\vec{1}, f) = 0$ (of course the domains of the sums are different.) By induction, we can assume that there exists a solution $B_g$ with $\nabla \cdot B_g = g$ on $T^-$ and we just want to extend $B_g$ to all of $T$. We define $A_f(a,b) = -f(a)$ and we just have to show that this definition satisfies $\nabla \cdot A_f(b) = f(b)$. But

$$\nabla \cdot A_f(b) = \sum_{(c,b) \in E} A_f(c,b) - \sum_{(b,c) \in E} A_f(b,c)$$

$$= A_f(a,b) + \sum_{(c,b) \in E^-} B_g(c,b) - \sum_{(b,c) \in E^-} B_g(b,c)$$

$$= -f(a) + \nabla \cdot B_g(b) = -f(a) + g(b) = -f(a) + f(b) + f(a) = f(b).$$

This shows that there is always a feasible solution $A_f$. Now we write $A = A_f + B$. Since $B$ is divergence free, it can be written as a linear combination of the functions in a cycle basis. Let $\{C_i \mid 1 \le i \le p\}$ be a cycle basis and write

$$B = \sum_i a_i C_i.$$

Then the energy $\xi(A_f + B) = \|A_f + B\|^2 = \|A_f\|^2 + 2(A_f, B) + (B, B)$ and the partial derivatives give

$$\frac{\partial \xi}{\partial a_i} = 2\big((B, C_i) + (A_f, C_i)\big) = \big(\nabla \times (A_f + B)\big)(C_i).$$

Thus at the extremal $\nabla \times A = \nabla \times (A_f + B) = 0$.

**Cycle updates**. The solution used in [3] to Equation 4 is based on cycle updates. We examine these updates in our language. We show that the update process never increases the energy $\xi$.

*Definition.* Let $W$ be a cycle basis. Given $A$ with $\nabla \cdot A = f$, consider the problem

$$\min_\alpha \|A + \alpha C'\|^2$$

for some fixed $C'$ in $W$. The solution is given by

$$\alpha^* = -\frac{(A, C')}{\|C'\|^2}.$$

We call $A + \alpha^* C'$ a *cycle update*.

*Lemma 3.* If $\nabla \cdot A = f$, then $\nabla \cdot (A + \alpha^* C') = f$ and

$$\xi(A + \alpha^* C') - \xi(A) = -\frac{(A, C')^2}{\|C'\|^2}.$$

*In particular, energy never increases in a cycle update.*

*Definition.* One useful type of cycle basis is associated with a collection of edges $T$ in $G$ that, considered as undirected, form a spanning tree for the graph, also considered as undirected. In this case each basis element is formed from a single edge in the graph but not in the tree together with the unique edges in the tree necessary to create a cycle. We call these cycles, *tree cycles,* and we call the cycle functions on these cycles, *tree cycle functions.* For each edge $e$ not in the tree, we denote the associated *tree cycle function* by $C_e$.

*Lemma 4. The tree cycle functions are a basis for the space of divergence free functions and the dimension is $p = m - n + 1$. This is also true when all entries are interpreted modulo a prime.*

*Proof.* Suppose that

$$\sum_{e \notin T} a_e C_e = 0.$$

Then only the tree cycle function $C_e$ has a non-zero entry at $e$ so $a_e = 0$. This shows the tree cycle functions are independent, both over the reals and the integers modulo a prime. The number of them is $m - n + 1$. In the $m$-dimensional space of functions on the edges of the graph, the indicator functions $T_e$ of the edges of the tree $T$ are independent and no linear combination of them can be divergence free. Thus the space of divergence free functions can have dimension no greater than $m - n + 1$ so the tree cycle functions are a basis.

*Lemma 5 (Duality.) Let $\nabla \cdot A = f$ on $G$ and $u$ be any function on the vertices. Then*

$$\|A\|^2 - (2(u, f) - (\nabla u, \nabla u)) = (A - \nabla u, A - \nabla u) \geq 0.$$

*Proof.* We use the standard argument:

$$0 \leq (A - \nabla u, A - \nabla u) = (A, A) - 2(\nabla u, A) + (\nabla u, \nabla u)$$

so

$$\|A\|^2 \geq 2(\nabla u, A) - (\nabla u, \nabla u) = 2(u, \nabla \cdot A) - (\nabla u, \nabla u).$$

*Definition. Let $\nabla \cdot A = f$ on $G$, $T$ be a spanning tree, $s$ be a leaf and $u$ be the induced function on the vertices.* Let $T'$ denote the set of edges not in $T$. Let

$$gap(A) = \|A\|^2 - (2(u,f) - (\nabla u, \nabla u)).$$

*Lemma 6. Let $\nabla \cdot A = f$ on $G$, $T$ be a spanning tree, $s$ be a leaf and $u$ be the induced function on the vertices. Then*

$$gap(A) = (A - \nabla u, A - \nabla u)_{T'}.$$

*Proof.* This follows directly from Lemma 5 and the fact that $A = \nabla u$ on $T$.

*Lemma 7 (tree cycle gap). Let $\nabla \cdot A = f$ on $G$, $T$ be a spanning tree, $s$ be a leaf and $u$ be the induced function on the vertices.*

*Then*

$$gap(A) = \sum_{e \in T'} (A, C_e)^2.$$

*Proof.* Let $C_e$ be a tree cycle function. Then since the curl of $\nabla u$ is zero

$$(A, C_e) = (A - \nabla u, C_e) = (A - \nabla u)(e).$$

Thus

$$(A - \nabla u, A - \nabla u)_{T'} = \sum_{e \in E \setminus T} (A, C_e)^2$$

which proves the lemma.

**General cycles**. Now we examine general bases with *pure* cycles, sets of edges that when direction is ignored form a cycle in the undirected graph. Let $\{D_i\}$ be a cycle basis that consists of pure cycles. Let $\nabla \cdot A = f$ on $G$, $T$ be a spanning tree, $s$ be a leaf and $u$ be the induced function on the vertices. Let $M$ be the incidence matrix of the edges relative to the (directed) cycles. Let $B = A - \nabla u$. Then we have seen that the error gap is

$$gap(A) = (B, B)_{T'}.$$

Since $\sum_i (B, D_i)^2 = (MB, MB)_{T'}$ and

$$(MB, MB) - (B, B) = (M^T MB, B) - (B, B) = ((M^T M - I)B, B),$$

if $M^T M - I$ is positive semi-definite, then

$$gap(A) = (B, B) \leq (MB, MB) = \sum_i (B, D_i)^2 = \sum_i (A, D_i)^2.$$

We state this as a theorem.

*Theorem 8 (cycle gap). Let $\{D_i\}$ be a cycle basis that consists of pure cycles. Let $\nabla \cdot A = f$ on $G$, $T$ be a spanning tree, $s$ be a leaf and $u$ be the induced function on the vertices. Let $M$ be the incidence matrix of the edges relative to the (directed) cycles. If $M^T M - I$ is positive semi-definite, then*

$$gap(A) \leq \sum_i (A, D_i)^2.$$

*Note.* Many cycle bases have this property. I don't know if all bases of pure cycles have this property.

**Expected convergence rate.** The iterative algorithm starts with a spanning tree $T$ and a current solution $A$ to $\nabla \cdot A = f$, picks a tree cycle function $C$ at random and updates to $A + \alpha^* C$ by a cycle update. Let $\tau = \sum_{e \in E \setminus T} \|C_e\|^2$. The probability that is used is

$$p(C) = \frac{1}{\tau} \|C\|^2.$$

*Lemma 9 (Expected progress). For a single cycle update, the expected value*

$$E(\xi(A + \alpha^* C) - \xi(A)) = -\frac{1}{\tau} gap(A).$$

*Proof.* We have by Lemma 3,

$$E(\xi(A + \alpha^* C) - \xi(A)) = -\sum_{e \in E \setminus T} p(C_e) \frac{(A, C_e)^2}{\|C_e\|^2} = -\frac{1}{\tau} \sum_{e \in E \setminus T} (A, C_e)^2.$$

*Theorem 10. The expected convergence rate of the energy $\xi(A)$ for the probabilistic algorithm is bounded by $1 - \frac{1}{\tau}$.*

*Proof.* Let $B$ a solution to Equation 4. Let $\nabla \cdot A = f$ on $G$, $T$ be a spanning tree, $s$ be a leaf and $u$ be the induced function on the vertices. Let

$$D(A) = \xi(A) - \xi(B).$$

By Lemma 6, $D(A) \leq gap(A)$. Thus using Lemma 9, we have

$$\frac{E(D(A+C_e))}{E(D(A))} = 1 - \frac{E(D(A)) - D(A+C_e)}{E(D(A))} = 1 - \frac{gap(A)}{\tau E(D(A))} \leq 1 - \frac{1}{\tau}.$$

*Note on low stretch trees.* There exists trees with $\tau = O(m \log n \log \log n)$. Utilizing one of these low stretch trees produces an algorithm that takes an expected $O(m \log n \log \log n)$ number of iterations to reduce the error in the energy by a fixed ratio. These trees are also used in [3] to show that the amount of work needed for a single cycle update is $O(\log n)$.

**Discussion of the probabilistic algorithm**. The point of the probability that is chosen is that longer cycles are updated more frequently than shorter cycles. If we imagine a traditional algorithm that just updates each cycle in some order in a single sweep and then repeats the sweep, each cycle would be updated the same number of times as any other. The convergence rate of the traditional algorithm would be much slower than the probabilistic algorithm because the longer cycles contribute more to the error than the shorter cycles. We could try to create an order for the cycles that includes the longer cycles proportionately more often than the shorter cycles but it is harder to see what that order might be. For example, just repeating a cycle a second time, does not reduce the error at all. If there was a logarithmic bound on the length of the longest cycle, a deterministic algorithm would produce an almost linear convergence scheme. However, there are practical graphs (for example, toroidal meshes) where this is not possible. If the number of long cycles in the cycle basis was bounded, we could produce a deterministic algorithm by updating all the long cycles after every cycle update. I don't know if bases like that always exist. For example, the toroidal meshes have this property; only two long cycles are required.